\documentclass[10pt,reqno]{amsart}
  \usepackage{geometry}
\usepackage{amsmath}
\usepackage{amsfonts}
\usepackage{amssymb}
\usepackage{graphicx}%
\providecommand{\U}[1]{\protect\rule{.1in}{.1in}}
\newtheorem{theorem}{Theorem}
\newtheorem{acknowledgement}[theorem]{Acknowledgement}

\newtheorem{corollary}[theorem]{Corollary}

\newtheorem{lemma}[theorem]{Lemma}

\newtheorem{proposition}[theorem]{Proposition}
\newtheorem{remark}[theorem]{Remark}

\begin{document}

\title{Orthogonality of the Dickson polynomials of the $(k+1)$-th kind }
\author{Diego Dominici}
\address{Johannes Kepler University Linz, Doktoratskolleg \textquotedblleft Computational Mathematics\textquotedblright, Altenberger Stra\ss e 69, 4040 Linz, Austria.}
\email{diego.dominici@dk-compmath.jku.at}
\address{Permanent address: Department of Mathematics, State University of New York at New Paltz, 1 Hawk Dr., New Paltz, NY 12561-2443, USA.}
\maketitle

\begin{abstract}
We study the Dickson polynomials of the $(k+1)$-th kind over the field of
complex numbers. We show that they are a family of co-recursive orthogonal
polynomials with respect to a quasi-definite moment functional $L_{k}.$ We
find an integral representation for $L_{k}$ and compute explicit expressions
for all of its moments.

\end{abstract}

\begin{quotation}
Dedicated to V. E. P., in gratitude for her continuous unbounded support
(without measure).
\end{quotation}


Classification codes: 11T06 (primary), 33C45, 26A42, 28C05 (secondary).

\section{Introduction}

Let $n\in\mathbb{N},$ $\mathbb{F}_{q}$ be a finite field and $a\in
\mathbb{F}_{q}.$ The Dickson polynomials $D_{n}(x;a),$ defined by
\cite[9.6.1]{MR3087321}
\[
D_{n}(x;a)=%
{\displaystyle\sum\limits_{j=0}^{\left\lfloor \frac{n}{2}\right\rfloor }}
\frac{n}{n-j}\binom{n-j}{j}\left(  -a\right)  ^{j}x^{n-2j},\quad
x\in\mathbb{F}_{q}%
\]
where introduced by Leonard Eugene Dickson (1874 -- 1954) in his 1896 Ph.D.
thesis \textquotedblleft The analytic representation of substitutions on a
power of a prime number of letters, with a discussion of the linear
group\textquotedblright\ \cite{MR2936778}, published in two parts in The
Annals of Mathematics \cite{MR1502214}, \cite{MR1502221}. The Dickson
polynomials are the unique monic polynomials satisfying the functional
equation \cite[9.6.3]{MR3087321}
\[
D_{n}\left(  y+\frac{a}{y};a\right)  =y^{n}+\left(  \frac{a}{y}\right)
^{n},\quad y\in\mathbb{F}_{q^{2}}.
\]
See \cite{MR1237403} for further algebraic and number theoretic properties of
the Dickson polynomials.

Let $\mathbb{N}_{0}$ denote the set $\mathbb{N}\cup\left\{  0\right\}
=0,1,2,\ldots.$ In \cite{MR2928474}, Wang and Yucas extended the Dickson
polynomials to a family depending on a new parameter $k\in\mathbb{N}_{0},$
which they called \emph{Dickson polynomials of the }$(k+1)$\emph{-th kind.
}They defined them by%

\begin{equation}
D_{n,k}(x;a)=%
{\displaystyle\sum\limits_{j=0}^{\left\lfloor \frac{n}{2}\right\rfloor }}
\frac{n-kj}{n-j}\binom{n-j}{j}\left(  -a\right)  ^{j}x^{n-2j},
\label{definition}%
\end{equation}
with initial values%
\begin{equation}
D_{0,k}(x;a)=2-k,\quad D_{1,k}(x;a)=x. \label{Initial}%
\end{equation}
They also showed that the polynomials $D_{n,k}(x;a)$ satisfy the fundamental
functional equation%
\begin{equation}
D_{n,k}(y+\frac{a}{y};a)=y^{n}+\left(  \frac{a}{y}\right)  ^{n}+k\frac
{ay^{n}-y^{2}\left(  \frac{a}{y}\right)  ^{n}}{y^{2}-a},\quad y\neq0.
\label{Functional}%
\end{equation}
Note that%
\[
\underset{y\rightarrow\pm\sqrt{a}}{\lim}\frac{ay^{n}-y^{2}\left(  \frac{a}%
{y}\right)  ^{n}}{y^{2}-a}=\left(  n-1\right)  \left(  \pm\sqrt{a}\right)
^{n},
\]
and therefore%
\[
D_{n,k}(y+\frac{a}{y};a)=\left[  2+\left(  n-1\right)  k\right]  \left(
\pm\sqrt{a}\right)  ^{n},\quad y=\pm\sqrt{a}.
\]

We clearly have%
\[
D_{n,0}(x;a)=D_{n}(x;a)\qquad\text{(Dickson polynomials)}%
\]
and \cite[9.6.1]{MR3087321}
\[
D_{n,1}(x;a)=E_{n}(x;a)\qquad\text{(Dickson polynomials of the second kind).}%
\]
In fact, since
\[
\frac{n-kj}{n-j}=k-\left(  k-1\right)  \frac{n}{n-j},
\]
we have \cite[2.1]{MR2928474}%
\[
D_{n,k}(x;a)=kE_{n}(x;a)-\left(  k-1\right)  D_{n}(x;a).
\]

The polynomials $D_{n,k}(x;a)$ also satisfy the fundamental recurrence (see
\cite[Remark 2.5]{MR2928474})%
\begin{equation}
D_{n+2,k}=xD_{n+1,k}-aD_{n,k},\quad n\in\mathbb{N}_{0}.\label{Reqn}%
\end{equation}
The first few Dickson polynomials of the $(k+1)$-th kind are%
\begin{align*}
D_{2,k}(x;a) &  =x^{2}+a\left(  k-2\right)  ,\\
D_{3,k}(x;a) &  =x^{3}+a\left(  k-3\right)  x,\\
D_{4,k}(x;a) &  =x^{4}+a\left(  k-4\right)  x^{2}+a^{2}\left(  2-k\right)  ,\\
D_{5,k}(x;a) &  =x^{5}+a\left(  k-5\right)  x^{3}+a^{2}\left(  5-2k\right)  x.
\end{align*}
They have zeros at%
\begin{align}
x &  =\pm\sqrt{a}\sqrt{2-k},\quad\text{if }n=2,\nonumber\\
x &  =0,\quad\pm\sqrt{a}\sqrt{3-k},\quad\text{if }n=3,\label{zeros}\\
x &  =\pm\frac{\sqrt{a}}{\sqrt{2}}\sqrt{4-k\pm\sqrt{\left(  k-2\right)
^{2}+4}},\quad\text{if }n=4,\nonumber\\
x &  =0,\quad\pm\frac{\sqrt{a}}{\sqrt{2}}\sqrt{5-k\pm\sqrt{\left(  k-1\right)
^{2}+4}},\quad\text{if }n=5,\nonumber
\end{align}
as can be verified using a mathematical symbolic computation program such as Mathematica.

\begin{remark}
Note that the polynomials $D_{n,k}(x;a)$ are monic for $n\geq1,$ but
$D_{0,k}(x;a)=1$ only for $k=1.$ 
\end{remark}

In this article, we study the polynomials $D_{n,k}(x;a)$ over the field of
complex numbers, with $a>0$ and $k\in\mathbb{R}.$ Our motivation is the
three-term recurrence relation (\ref{Reqn}), which suggests that the Dickson
polynomials of the $(k+1)$-th kind form a family of orthogonal polynomials
with respect to some linear functional $L_{k}$. However, from (\ref{zeros}) we
see that for $k>2$ the polynomials $D_{n,k}(x;a)$ may have a pair of purely
imaginary roots. Also, the polynomials $D_{3,3}(x;a)$ and $D_{5,\frac{5}{2}%
}(x;a)$ have a triple zero at $x=0.$ This implies that the linear functional
$L_{k}$ is quasi-definite \cite[Theorem 2.4.3]{MR2703093}, \cite[Theorem
1]{MR1112321}.

The article is organized as follows: in Section \ref{Section1}, we derive some
of the main properties of the Dickson polynomials of the $(k+1)$-th kind,
including different expressions, a hypergeometric representation, differential
equations, and a generating function.

In Section \ref{Section2}, we present some basic results from the theory of
orthogonal polynomials that we will need to find the linear functional $L_{k}%
$. We define the co-recursive polynomials associated with a given family of
orthogonal polynomials, and list some of their main properties. We also show
that a family polynomials related to $D_{n,k}(x;a)$ are co-recursive
polynomials associated with the Chebyshev polynomials of the second kind.

In Section \ref{Section3} we apply the results of the previous sections to the
Dickson polynomials of the $(k+1)$-th kind and obtain a representation for
their moment functional $L_{k}.$ We also find explicit expressions for the
moments of $L_{k}.$

Finally, in Section \ref{Section4} we summarize our results. In our hope that
the results would be of interest to researchers outside the field of
orthogonal polynomials and special functions, we have made the paper as
self-contained as possible.

\section{Properties of Dickson polynomials \label{Section1}}

We begin by checking the initial polynomial $D_{0,k}(x;a).$ Since it is not
clear from the definition (\ref{definition}) that $D_{0,k}(x;a)=2-k,$ we
consider even and odd degrees and obtain the following result.

\begin{proposition}
The even and odd Dickson polynomials of the $(k+1)$-th kind are given by%
\begin{equation}
D_{2n,k}(x;a)=\left(  2-k\right)  \left(  -a\right)  ^{n}+%
{\displaystyle\sum\limits_{j=1}^{n}}
\frac{\left(  2-k\right)  n+kj}{j+n}\binom{n+j}{2j}\left(  -a\right)
^{n-j}x^{2j} \label{even}%
\end{equation}
and%
\begin{equation}
D_{2n+1,k}(x;a)=x%
{\displaystyle\sum\limits_{j=0}^{n}}
\frac{\left(  2-k\right)  n+kj+1}{j+n+1}\binom{n+j+1}{2j+1}\left(  -a\right)
^{n-j}x^{2j}. \label{odd}%
\end{equation}

\end{proposition}

\begin{proof}
From (\ref{definition}), we have%
\[
D_{2n,k}(x;a)=%
{\displaystyle\sum\limits_{j=0}^{n}}
\frac{2n-kj}{2n-j}\binom{2n-j}{j}\left(  -a\right)  ^{j}x^{2n-2j},
\]
and switching the index to $i=n-j,$ we get
\[
D_{2n,k}(x;a)=%
{\displaystyle\sum\limits_{i=0}^{n}}
\frac{\left(  2-k\right)  n+ki}{i+n}\binom{n+i}{n-i}\left(  -a\right)
^{n-i}x^{2i},
\]
and (\ref{even}) follows after using the symmetry of the binomial coefficients%
\[
\binom{n}{k}=\binom{n}{n-k}.
\]

A similar calculation gives (\ref{odd}).
\end{proof}

Next, we will find a representation for $D_{n,k}(x;a)$ in terms of the
generalized hypergeometric function
\[
_{p}F_{q}\left(
\begin{array}
[c]{c}%
a_{1},\ldots,a_{p}\\
b_{1},\ldots,b_{q}%
\end{array}
;x\right)  =%
{\displaystyle\sum\limits_{j=0}^{\infty}}
\frac{\left(  a_{1}\right)  _{j}\cdots\left(  a_{p}\right)  _{j}}{\left(
b_{1}\right)  _{j}\cdots\left(  b_{q}\right)  _{j}}\frac{x^{j}}{j!},
\]
where $\left(  u\right)  _{j}$ denotes the Pochhammer symbol (also called
shifted or rising factorial) defined by \cite[5.2.4]{MR2723248}
\begin{align*}
\left(  a\right)  _{0}  &  =1\\
\left(  a\right)  _{j}  &  =a\left(  a+1\right)  \cdots\left(  a+j-1\right)
,\quad j\in\mathbb{N}.
\end{align*}

\begin{proposition}
The Dickson polynomials of the $(k+1)$-th kind admit the hypergeometric
representation%
\begin{equation}
D_{n,k}(x;a)=x^{n}\ _{3}F_{2}\left(
\begin{array}
[c]{c}%
-\frac{n-1}{2},-\frac{n}{2},1-\frac{n}{k}\\
1-n,-\frac{n}{k}%
\end{array}
;\frac{4a}{x^{2}}\right)  ,\quad k\neq0. \label{hyper}%
\end{equation}
We also have%
\[
D_{n,0}(x;a)=x^{n}\ _{2}F_{1}\left(
\begin{array}
[c]{c}%
-\frac{n-1}{2},-\frac{n}{2}\\
1-n
\end{array}
;\frac{4a}{x^{2}}\right)  .
\]

\end{proposition}

\begin{proof}
Let
\[
D_{n,k}(x;a)=x^{n}%
{\displaystyle\sum\limits_{j=0}^{\left\lfloor \frac{n}{2}\right\rfloor }}
c_{j},
\]
with
\begin{equation}
c_{j}=\frac{n-kj}{n-j}\binom{n-j}{j}\left(  -a\right)  ^{j}x^{-2j}. \label{cj}%
\end{equation}
We have $c_{0}=1$ and $c_{j}=0$ for $j>\frac{n}{2}.$ Using (\ref{cj}), we get%
\begin{equation}
\frac{c_{j+1}}{c_{j}}=\frac{\left(  2j-n+1\right)  \left(  2j-n\right)
\left(  jk+k-n\right)  }{\left(  j-n+1\right)  \left(  jk-n\right)  \left(
j+1\right)  }\frac{a}{x^{2}}. \label{ratio}%
\end{equation}

Let $k\neq0.$ Then,
\[
\frac{c_{j+1}}{c_{j}}=\frac{\left(  j-\frac{n}{2}+\frac{1}{2}\right)  \left(
j-\frac{n}{2}\right)  \left(  j+1-\frac{n}{k}\right)  }{\left(  j-n+1\right)
\left(  j-\frac{n}{k}\right)  \left(  j+1\right)  }\frac{4a}{x^{2}},
\]
and it follows that%
\[
c_{j}=\frac{\left(  -\frac{n-1}{2}\right)  _{j}\left(  -\frac{n}{2}\right)
_{j}\left(  1-\frac{n}{k}\right)  _{j}}{\left(  1-n\right)  _{j}\left(
-\frac{n}{k}\right)  _{j}}\frac{1}{j!}\left(  \frac{4a}{x^{2}}\right)  ^{j}.
\]
Thus,
\[
D_{n,k}(x;a)=x^{n}%
{\displaystyle\sum\limits_{j=0}^{\infty}}
\frac{\left(  -\frac{n-1}{2}\right)  _{j}\left(  -\frac{n}{2}\right)
_{j}\left(  1-\frac{n}{k}\right)  _{j}}{\left(  1-n\right)  _{j}\left(
-\frac{n}{k}\right)  _{j}}\frac{1}{j!}\left(  \frac{4a}{x^{2}}\right)  ^{j}.
\]

If $k=0,$ we see from (\ref{ratio}) that%
\[
\frac{c_{j+1}}{c_{j}}=\frac{4\left(  j+\frac{1-n}{2}\right)  \left(
j-\frac{n}{2}\right)  }{\left(  j-n+1\right)  \left(  j+1\right)  }\frac
{a}{x^{2}},
\]
and therefore%
\[
c_{j}=\frac{\left(  -\frac{n-1}{2}\right)  _{j}\left(  -\frac{n}{2}\right)
_{j}}{\left(  1-n\right)  _{j}}\frac{1}{j!}\left(  \frac{4a}{x^{2}}\right)
^{j}.
\]
Hence,%
\[
D_{n,0}(x;a)=x^{n}%
{\displaystyle\sum\limits_{j=0}^{\infty}}
\frac{\left(  -\frac{n-1}{2}\right)  _{j}\left(  -\frac{n}{2}\right)  _{j}%
}{\left(  1-n\right)  _{j}}\frac{1}{j!}\left(  \frac{4a}{x^{2}}\right)  ^{j}.
\]

\end{proof}

\begin{remark}
A representation of $D_{n,2}(x;a)$ in terms of associated Legendre functions
of the first and second kinds \cite[14.3]{MR2723248} was given in "A
representation of the Dickson polynomials of the third kind by Legendre
functions" by Neranga Fernando and Solomon Manukure (arXiv:1604.04682).
\end{remark}

\begin{proposition}
For $n\in\mathbb{N}_{0},$ the Dickson polynomials of the $(k+1)$-th kind
satisfy the following relations:%
\[
D_{n,k+2}-2D_{n,k+1}+D_{n,k}=0,
\]%
\begin{gather}
-(x^{2}-4a)\left[  (k-1)nx^{2}+a(k-2)(kn-2n-k)\right]  D_{n,k}^{\prime\prime
}\nonumber\\
+x\left[  (k-1)nx^{2}+a\left(  6k+4n+3k^{2}n-4kn-3k^{2}\right)  \right]
D_{n,k}^{\prime}\label{ODE}\\
+\left[  (k-1)n^{3}x^{2}+an\left(  -k-2n+kn\right)  \left(  -2k-2n+kn\right)
\right]  D_{n,k}=0,\nonumber
\end{gather}%
\begin{gather*}
(x^{2}-4a)^{2}D_{n,k}^{\left(  iv\right)  }+10x(x^{2}-4a)D_{n,k}^{\prime
\prime\prime}+\left[  (23-2n^{2})x^{2}+8a(n^{2}-4)\right]  D_{n,k}%
^{\prime\prime}\\
-3(2n^{2}-3)xD_{n,k}^{\prime}+n^{2}(n^{2}-4)D_{n,k}=0,
\end{gather*}
and%
\[
(x^{2}-4a)D_{n,k}^{\prime\prime}-4nD_{n+1,k}D_{n,k}^{\prime}+\left(
2n+3\right)  xD_{n,k}^{\prime}+n\left(  n+2\right)  D_{n,k}=0.
\]

\end{proposition}

\begin{proof}
All the identities can be automatically found and proved using the
hypergeometric representation (\ref{hyper}) and the Mathematica package
HolonomicFunctions \cite{MR3624284}.
\end{proof}

\begin{remark}
The differential equation (\ref{ODE}) already appeared in \cite[Lemma
2.7]{MR2928474}.
\end{remark}

We can use the recurrence relation (\ref{Reqn}) to obtain a different
representation for the polynomials $D_{n,k}(x;a).$

\begin{proposition}
For $x\neq\pm2\sqrt{a},$ the Dickson polynomials of the $(k+1)$-th kind are
given by%
\begin{equation}
D_{n,k}(x;a)=\left(  1+k\frac{x-\Delta}{2\Delta}\right)  \left(
\frac{x+\Delta}{2}\right)  ^{n}+\left(  1-k\frac{x+\Delta}{2\Delta}\right)
\left(  \frac{x-\Delta}{2}\right)  ^{n}, \label{Dnk}%
\end{equation}
where%
\[
\Delta=\sqrt{x^{2}-4a}.
\]

We also have,%
\begin{equation}
D_{n,k}(\pm2\sqrt{a};a)=\left(  kn+2-k\right)  \left(  \pm\sqrt{a}\right)
^{n}. \label{Dnkroot}%
\end{equation}

\end{proposition}

\begin{proof}
Let us assume that we can write%
\begin{equation}
D_{n,k}(x;a)=R^{n}, \label{R}%
\end{equation}
for some function $R(x,k,a).$ Using (\ref{R}) in the recurrence (\ref{Reqn}),
we obtain
\[
R^{2}-xR+a=0,
\]
and therefore%
\[
R_{\pm}=\frac{x\pm\Delta}{2}%
\]
with%
\[
\Delta=\sqrt{x^{2}-4a}.
\]
It follows that the general solution of (\ref{Reqn}) is given by%
\begin{equation}
D_{n,k}(x;a)=C_{1}\left(  x;a,k\right)  \left(  \frac{x+\Delta}{2}\right)
^{n}+C_{2}\left(  x;a,k\right)  \left(  \frac{x-\Delta}{2}\right)  ^{n}.
\label{General}%
\end{equation}

Using the initial conditions (\ref{Initial}) in (\ref{General}), we get%
\begin{align*}
C_{1}\left(  x;a,k\right)  +C_{2}\left(  x;a,k\right)   &  =2-k,\\
C_{1}\left(  x;a,k\right)  \left(  \frac{x+\Delta}{2}\right)  +C_{2}\left(
x;a,k\right)  \left(  \frac{x-\Delta}{2}\right)   &  =x.
\end{align*}
Thus, assuming that $x\neq\pm2\sqrt{a},$
\[
C_{1}\left(  x;a,k\right)  =1+k\frac{x-\Delta}{2\Delta},\quad C_{2}\left(
x;a,k\right)  =1-k\frac{x+\Delta}{2\Delta}.
\]

To verify (\ref{Dnkroot}), we replace it in the recurrence (\ref{Reqn}), and
obtain%
\begin{align*}
&  \left(  kn+2+k\right)  \left(  \pm\sqrt{a}\right)  ^{n+2}-\left(  \pm
2\sqrt{a}\right)  \left(  kn+2\right)  \left(  \pm\sqrt{a}\right)
^{n+1}+a\left(  kn+2-k\right)  \left(  \pm\sqrt{a}\right)  ^{n}\\
&  =a\left(  \pm\sqrt{a}\right)  ^{n}\left[  \left(  kn+2+k\right)  -2\left(
kn+2\right)  +\left(  kn+2-k\right)  \right]  =0.
\end{align*}

\end{proof}

Using (\ref{Dnk}), we can obtain a generating function for the polynomials
$D_{n,k}(x;a).$

\begin{proposition}
The ordinary generating function of the polynomials $D_{n,k}(x;a)$ is given
by
\begin{equation}
G(z;x,k,a)=%
{\displaystyle\sum\limits_{n=0}^{\infty}}
D_{n,k}(x;a)z^{n}=\frac{2-k+\left(  k-1\right)  xz}{az^{2}-xz+1}.
\label{GenFun}%
\end{equation}
$\allowbreak$
\end{proposition}

\begin{proof}
From (\ref{Dnk}), we have (as formal power series)%
\begin{align*}
G(z;x,k,a)  &  =%
{\displaystyle\sum\limits_{n=0}^{\infty}}
D_{n,k}(x;a)z^{n}\\
&  =\left(  1+k\frac{x-\Delta}{2\Delta}\right)
{\displaystyle\sum\limits_{n=0}^{\infty}}
\left(  z\frac{x+\Delta}{2}\right)  ^{n}+\left(  1-k\frac{x+\Delta}{2\Delta
}\right)
{\displaystyle\sum\limits_{n=0}^{\infty}}
\left(  z\frac{x-\Delta}{2}\right)  ^{n}\\
&  =\left(  1+k\frac{x-\Delta}{2\Delta}\right)  \frac{1}{1-\left(
z\frac{x+\Delta}{2}\right)  }+\left(  1-k\frac{x+\Delta}{2\Delta}\right)
\frac{1}{1-\left(  z\frac{x-\Delta}{2}\right)  }.
\end{align*}
Thus,%
\[
G(z;x,k,a)=4\frac{2-k+\left(  k-1\right)  xz}{\left(  x^{2}-\Delta^{2}\right)
z^{2}-4xz+4}%
\]
$\allowbreak$and the result follows since
\[
\Delta^{2}=x^{2}-4a.
\]

\end{proof}

\begin{remark}
The same generating function was obtained in \cite[Lemma 2.6]{MR2928474} using
the recurrence (\ref{Reqn}).
\end{remark}

\section{Orthogonal polynomials \label{Section2}}

Let $\left\{  \mu_{n}\right\}  $ be a sequence of complex numbers and
$\mathcal{L}:\mathbb{C}\left[  x\right]  \rightarrow\mathbb{C}$ be a linear
functional defined by%
\[
\mathcal{L}\left[  x^{n}\right]  =\mu_{n},\quad n\in\mathbb{N}_{0}.
\]
Then, $\mathcal{L}$ is called the \emph{moment functional} determined by the
moment sequence $\left\{  \mu_{n}\right\}  .$ The number $\mu_{n}$ is called
the \emph{moment} of order $n$.

A moment functional $\mathcal{L}$ is called \emph{positive-definite} if
$\mathcal{L}\left[  r(x)\right]  >0$ for every polynomial $r(x)$ that is not
identically zero and is non-negative for all real $x$. Otherwise,
$\mathcal{L}$ is called \emph{quasi-definite}.

A sequence $\left\{  P_{n}\right\}  \subset\mathbb{C}\left[  x\right]  ,$ with
$\deg\left(  P_{n}\right)  =n$ is called an \emph{orthogonal polynomial
sequence} with respect to\ $\mathcal{L}$ provided that \cite{MR0481884}%
\[
\mathcal{L}\left[  P_{n}P_{m}\right]  =h_{n}\delta_{n,m},\quad n,m\in
\mathbb{N}_{0},
\]
where $h_{n}\neq0$ and $\delta_{n,m}$ is Kronecker's delta.

One of the fundamental properties of orthogonal polynomials is that they
satisfy a three-term recurrence relation.

\begin{theorem}
\label{Theorem1}Let $\mathcal{L}$ be a moment functional and let $\left\{
P_{n}\right\}  $ be the sequence of monic orthogonal polynomials associated
with it. Then, there exist $\beta_{n}\in\mathbb{C}$ and $\gamma_{n}%
\in\mathbb{C}\setminus\left\{  0\right\}  $ such that the polynomials
$P_{n}\left(  x\right)  $ are a solution of the 3-term recurrence relation
\begin{equation}
P_{n+1}=\left(  x-\beta_{n}\right)  P_{n}-\gamma_{n}P_{n-1},\quad
n\in\mathbb{N}, \label{3-term}%
\end{equation}
with initial conditions%
\begin{equation}
P_{0}\left(  x\right)  =1,\quad P_{1}\left(  x\right)  =x-\beta_{0}.
\label{initialP}%
\end{equation}

\end{theorem}

\begin{proof}
See \cite[Theorem 4.1]{MR0481884}.
\end{proof}

A linearly independent solution of (\ref{3-term}) with initial conditions%
\begin{equation}
P_{0}^{\ast}\left(  x\right)  =0,\quad P_{1}^{\ast}\left(  x\right)  =1,
\label{initialPa}%
\end{equation}
is given by the so-called \emph{associated orthogonal polynomials}
$P_{n}^{\ast}\left(  x\right)  $ \cite[4.3]{MR0481884}. Note that $\deg
P_{n}^{\ast}\left(  x\right)  =n-1.$

The converse of Theorem \ref{Theorem1} is known as Favard's Theorem.

\begin{theorem}
Let $\left\{  P_{n}\right\}  $ be a sequence of polynomials satisfying the
3-term recurrence relation (\ref{3-term}) with $\beta_{n}\in\mathbb{C}$ and
$\gamma_{n}\in\mathbb{C}\setminus\left\{  0\right\}  .$ Then, there exists a
\underline{unique} linear functional $\mathcal{L}$ such that
\[
\mathcal{L}\left[  P_{0}\right]  =1,\quad\mathcal{L}\left[  P_{n}P_{m}\right]
=h_{n}\delta_{n,m},
\]
with
\begin{equation}
h_{0}=P_{0},\quad h_{1}=\gamma_{1},\quad h_{n}=\gamma_{n}h_{n-1},\quad
n=2,3,\ldots. \label{hn}%
\end{equation}

\end{theorem}

\begin{proof}
See \cite[Theorem 4.4]{MR0481884}.
\end{proof}

\begin{remark}
\label{RemarkL}It follows from (\ref{Initial}) and (\ref{Reqn}) that (at least
for $k\neq2)$ $\left\{  D_{n,k}\right\}  $ is a sequence of monic (for
$n\geq1)$ orthogonal polynomials with respect to a moment functional $L_{k}$
\footnote{In the remainder of the paper $L_{k}$ will denote the moment
functional associated with the polynomials $D_{n,k}(x;a).$} satisfying%
\begin{equation}
L_{k}\left[  D_{n,k}D_{m,k}\right]  =h_{n}\left(  k\right)  \delta
_{n,m},\label{L1}%
\end{equation}
with
\begin{equation}
h_{0}\left(  k\right)  =2-k,\quad h_{n}\left(  k\right)  =a^{n},\quad
n\in\mathbb{N}.\label{hD}%
\end{equation}

In Section \ref{Section3} we will find a representation for the moment
functional $L_{k}.$
\end{remark}

\begin{proposition}
\label{Prop1}Let $\mathcal{L}$ be a moment functional and $\left\{
P_{n}\right\}  $ be the sequence of monic orthogonal polynomials associated
with it. Then, the following are equivalent:

(a) All the moments of odd order are zero,
\[
\mathcal{L}\left[  x^{2n+1}\right]  =0,\quad n\in\mathbb{N}_{0}.
\]

(b) The polynomials $P_{n}\left(  x\right)  $ satisfy
\[
P_{n}\left(  -x\right)  =\left(  -1\right)  ^{n}P_{n}\left(  x\right)  ,\quad
n\in\mathbb{N}_{0}.
\]

\end{proposition}

\begin{proof}
See \cite[Theorem 4.3]{MR0481884}.
\end{proof}

\begin{proposition}
Let $k\neq2$ and $\mu_{n}\left(  k\right)  $ denote the moments of the linear
functional defined by (\ref{L1}). Then, we have%
\begin{equation}
\mu_{0}\left(  k\right)  =\frac{1}{2-k}, \label{mu0}%
\end{equation}%
\begin{equation}
\mu_{2n+1}\left(  k\right)  =0,\quad n\in\mathbb{N}_{0}, \label{oddmu}%
\end{equation}
and%
\[
\mu_{2n}=-%
{\displaystyle\sum\limits_{j=0}^{n-1}}
\frac{\left(  2-k\right)  n+kj}{j+n}\binom{n+j}{2j}\left(  -a\right)
^{n-j}\mu_{2j},\quad n\in\mathbb{N}.
\]

The first few nonzero moments are%
\begin{align*}
\mu_{2}\left(  k\right)   &  =a,\quad\mu_{4}\left(  k\right)  =-a^{2}\left(
k-3\right)  ,\quad\mu_{6}\left(  k\right)  =a^{3}(k^{2}-6k+10),\\
\mu_{8}\left(  k\right)   &  =-a^{4}(k^{3}-9k^{2}+29k-35).
\end{align*}

\end{proposition}

\begin{proof}
From (\ref{hD}), we see that%
\[
2-k=h_{0}=L_{k}\left[  D_{0,k}^{2}\right]  =D_{0,k}^{2}L_{k}\left[  1\right]
=\left(  2-k\right)  ^{2}\mu_{0},
\]
from which (\ref{mu0}) follows.

Using (\ref{definition}), it is clear that%
\begin{equation}
D_{n,k}\left(  -x;a\right)  =\left(  -1\right)  ^{n}D_{n,k}\left(  x;a\right)
, \label{oddeven}%
\end{equation}
and Proposition \ref{Prop1} gives
\[
\mu_{2n+1}\left(  k\right)  =0,\quad n=0,1,\ldots.
\]

From (\ref{even}) we have%
\[
D_{2n,k}(x;a)=%
{\displaystyle\sum\limits_{j=0}^{n}}
\frac{\left(  2-k\right)  n+kj}{j+n}\binom{n+j}{2j}\left(  -a\right)
^{n-j}x^{2j},
\]
and therefore%
\begin{align*}
0  &  =L_{k}\left[  D_{2n,k}\right]  =%
{\displaystyle\sum\limits_{j=0}^{n}}
\frac{\left(  2-k\right)  n+kj}{j+n}\binom{n+j}{2j}\left(  -a\right)
^{n-j}\mu_{2j}\left(  k\right) \\
&  =%
{\displaystyle\sum\limits_{j=0}^{n-1}}
\frac{\left(  2-k\right)  n+kj}{j+n}\binom{n+j}{2j}\left(  -a\right)
^{n-j}\mu_{2j}\left(  k\right)  +\mu_{2n}\left(  k\right)  .
\end{align*}

\end{proof}

\begin{remark}
In Section \ref{Section3} we will find a closed-form expression for $\mu
_{2n}\left(  k\right)  .$
\end{remark}

The task of finding an explicit integral representation for the functional
$\mathcal{L}$ is called a \emph{moment problem} \cite{MR0184042}%
,\cite{MR0458081},\cite{MR0008438}. A moment functional $\mathcal{L}$ is
called \emph{determinate} if there exists a unique (up to an additive
constant) distribution $\psi\left(  x\right)  $ such that%
\begin{equation}
\mathcal{L}\left[  x^{n}\right]  =%
{\displaystyle\int\limits_{\Lambda}}
x^{n}d\psi\left(  x\right)  , \label{mu}%
\end{equation}
where the set $\Lambda$ is called the support of the distribution $\psi.$
Otherwise, $\mathcal{L}$ is called \emph{indeterminate} \cite{MR1285262},
\cite{MR2242752}.

A criteria to decide if the moment functional $\mathcal{L}$ is determinate is
due to Torsten Carleman \cite[P 59]{MR0008438}: If $\gamma_{n}>0$ and
\[%
{\displaystyle\sum\limits_{n=1}^{\infty}}
\frac{1}{\sqrt{\gamma_{n}}}=\infty,
\]
then $\mathcal{L}$ is determinate. Since for the Dickson polynomials we have
$\gamma_{n}=a,$ it follows that the moment problem is determinate.

One method to find a distribution function satisfying (\ref{mu}) is given by
Markov's theorem.

\begin{theorem}
Let the moment functional $\mathcal{L},$ supported on the set $\Lambda
\subset\mathbb{C},$ be determinate, $\left\{  P_{n}\right\}  $ be the monic
orthogonal polynomials with respect to $\mathcal{L}$ and $\left\{  P_{n}%
^{\ast}\right\}  $ be the associated polynomials. Then,
\begin{equation}
\underset{n\rightarrow\infty}{\lim}\mathcal{L}\left[  1\right]  \frac
{P_{n}^{\ast}\left(  z\right)  }{P_{n}(z)}=\mathcal{L}\left[  \frac{1}%
{z-x}\right]  ,\quad z\notin\Lambda, \label{markov}%
\end{equation}
where the convergence is uniform on compact subsets of $\mathbb{C}%
\setminus\Lambda.$
\end{theorem}

\begin{proof}
See \cite{MR1136928}.
\end{proof}

The function
\begin{equation}
\mathcal{S}(z)=\mathcal{L}\left[  \frac{1}{z-x}\right]  , \label{Stieltjes}%
\end{equation}
is called the \emph{Stieltjes transform} of $\mathcal{L}$ \cite{MR903848} and
it has the asymptotic behavior \cite[Section 12.9]{MR0453984}
\begin{equation}
\mathcal{S}(z)\sim\frac{\mu_{0}}{z}+\frac{\mu_{1}}{z^{2}}+\frac{\mu_{2}}%
{z^{3}}+\cdots,\quad z\rightarrow\infty. \label{asympt}%
\end{equation}

\subsection{Co-recursive polynomials}

Let $\left\{  q_{n}\right\}  $ be a sequence of polynomials satisfying
(\ref{3-term}) with initial conditions%
\begin{equation}
q_{0}\left(  x\right)  =u,\quad q_{1}\left(  x\right)  =x-v. \label{Initialq}%
\end{equation}
The polynomials $q_{n}\left(  x\right)  $ are called \emph{co-recursive }with
parameters $\left(  u,v\right)  $. They were introduced by T. Chihara in
\cite{MR0092015}, where he considered the case $u=1.$ See also
\cite{MR2702951}, \cite{MR1266575}, \cite{MR0420120}, and \cite{MR972579}.

Note that when $u=1$ we could consider the polynomials $q_{n}\left(  x\right)
$ to be solutions of the perturbed 3-term recurrence relation%
\[
xq_{n}=q_{n+1}+\widetilde{\beta}_{n}q_{n}+\gamma_{n}q_{n-1},\quad
n\in\mathbb{N}_{0}%
\]
with initial conditions%
\[
q_{-1}\left(  x\right)  =0,\quad q_{0}\left(  x\right)  =1,
\]
where%
\[
\widetilde{\beta}_{n}=\beta_{n}+\left(  v-\beta_{0}\right)  \delta_{n,0}.
\]
Orthogonal polynomials that are solutions of 3-term recurrence relations with
finite perturbations of the recurrence coefficients are called
\emph{co-polynomials}. See \cite{MR3318210}, \cite{MR2005575}, and
\cite{MR1062324}.

Let $\left\{  P_{n}\right\}  $ and $\left\{  P_{n}^{\ast}\right\}  $ be the
linearly independent solutions of (\ref{3-term}) with initial conditions
(\ref{initialP}) and (\ref{initialPa}). Then, we can represent $q_{n}\left(
x\right)  $ as a linear combination%
\[
q_{n}\left(  x\right)  =AP_{n}\left(  x\right)  +BP_{n}^{\ast}\left(
x\right)  .
\]
Using (\ref{Initialq}), we see that%
\[
A=u,\quad B=\left(  1-u\right)  x+u\beta_{0}-v,
\]
and therefore%
\begin{equation}
q_{n}\left(  x\right)  =uP_{n}\left(  x\right)  +\left[  \left(  1-u\right)
x+u\beta_{0}-v\right]  P_{n}^{\ast}\left(  x\right)  ,\quad n\in\mathbb{N}%
_{0}. \label{q1}%
\end{equation}
Linear combinations of orthogonal polynomials have been studied by many
authors. See \cite{MR2010273}, \cite{MR1865881}, \cite{MR1414283},
\cite{MR1449688}, and \cite{MR3180027}.

Note that the associated polynomials for both sequences are the same, i.e.,
$q_{n}^{\ast}\left(  x\right)  =P_{n}^{\ast}\left(  x\right)  .$ Thus, we have%
\[
\frac{q_{n}\left(  z\right)  }{q_{n}^{\ast}\left(  z\right)  }=u\frac
{P_{n}\left(  z\right)  }{P_{n}^{\ast}\left(  z\right)  }+\left(  1-u\right)
z+u\beta_{0}-v,
\]
and assuming that the moment problem is determined, we can use Markov's
theorem and obtain%
\[
\frac{L_{q}\left[  1\right]  }{S_{q}\left(  z\right)  }=\frac{uL_{P}\left[
1\right]  }{S_{P}\left(  z\right)  }+\left(  1-u\right)  z+u\beta_{0}-v,
\]
where $S_{q}\left(  z\right)  ,$ $S_{p}\left(  z\right)  $ and $L_{q},$
$L_{p}$ denote the Stieltjes transforms and linear functionals associated with
the sequences $\left\{  q_{n}\right\}  $ and $\left\{  P_{n}\right\}  ,$
respectively. Solving for $S_{q}\left(  z\right)  ,$ we obtain%
\begin{equation}
S_{q}\left(  z\right)  =\frac{L_{q}\left[  1\right]  S_{P}\left(  z\right)
}{\left[  \left(  1-u\right)  z+u\beta_{0}-v\right]  S_{P}\left(  z\right)
+uL_{P}\left[  1\right]  }.\label{Sq}%
\end{equation}
The rational transformation (\ref{Sq}) is a particular case of the general
\emph{spectral transformations} of the Stieltjes transform \cite{MR1482157}%
\[
\widetilde{S}\left(  z\right)  =\frac{A\left(  z\right)  S\left(  z\right)
+B\left(  z\right)  }{C\left(  z\right)  S\left(  z\right)  +D\left(
z\right)  },
\]
where $A\left(  z\right)  ,B\left(  z\right)  ,C\left(  z\right)  ,D\left(
z\right)  $ are polynomials. See also \cite{MR2055354}, \cite{MR3008849},
\cite{MR1140648}, \cite{MR1934914}, and \cite{MR1199259}.

Suppose that $L_{P}\left[  1\right]  =1,$ and let's denote the polynomial
$\left(  1-u\right)  z+u\beta_{0}-v$ appearing in (\ref{Sq}) by $Q\left(
z\right)  .$ If the function $S_{P}\left(  z\right)  $ is a solution of the
quadratic equation%
\begin{equation}
AQS_{P}^{2}+\left(  BQ+uA-C\right)  S_{P}+uB=0, \label{quadratic}%
\end{equation}
for some polynomials $A\left(  z\right)  ,B\left(  z\right)  ,C\left(
z\right)  ,$ then we have%
\[
\frac{S_{q}\left(  z\right)  }{L_{q}\left[  1\right]  }=\frac{S_{P}\left(
z\right)  }{Q\left(  z\right)  S_{P}\left(  z\right)  +u}=\frac{A\left(
z\right)  S_{P}\left(  z\right)  +B\left(  z\right)  }{C\left(  z\right)  }.
\]
Stieltjes transforms satisfying quadratic equations with polynomial
coefficients are called \emph{second degree forms}. See \cite{MR2879183},
\cite{MR1755691}, \cite{MR1314902}, and \cite{MR2837719}.

The linear transformation%
\[
\frac{S_{q}\left(  z\right)  }{L_{q}\left[  1\right]  }=\frac{A\left(
z\right)  S_{P}\left(  z\right)  +B\left(  z\right)  }{C\left(  z\right)  }%
\]
can be written as a composition of three basic transformations:

\begin{enumerate}
\item The \emph{Uvarov transformation}, \cite{MR1199643}, \cite{MR0262764},
defined by%
\[
\frac{L_{U}\left[  r\right]  }{\lambda_{U}}=L_{P}\left[  r\right]  +Mr\left(
\omega\right)  ,\quad r\in\mathbb{C}\left[  x\right]  ,
\]
or by%
\[
\frac{S_{U}\left(  z\right)  }{\lambda_{U}}=S\left(  z\right)  +\frac
{M}{z-\omega},
\]
where $M+L_{P}\left[  1\right]  \neq0$ and%
\[
\lambda_{U}=\frac{L_{U}\left[  1\right]  }{L_{P}\left[  1\right]  +M}.
\]

\item The \emph{Christoffel transformation}, \cite{MR2324861},
\cite{MR1579059} defined by%
\[
\frac{L_{C}\left[  r\right]  }{\lambda_{C}}=L_{P}\left[  \left(
x-\omega\right)  r\left(  x\right)  \right]  ,\quad r\in\mathbb{C}\left[
x\right]  ,
\]
or by%
\[
\frac{S_{C}\left(  z\right)  }{\lambda_{C}}=\left(  z-\omega\right)  S\left(
z\right)  -L_{P}\left[  1\right]  ,
\]
where $L_{P}\left[  x-\omega\right]  \neq0$ and%
\[
\lambda_{C}=\frac{L_{C}\left[  1\right]  }{L_{P}\left[  x-\omega\right]  }.
\]

\item The \emph{Geronimus transformation }\cite{MR3208415}, \cite{MR3264577},
\cite{MR0004339}, \cite{MR0004340}, defined by%
\[
\frac{L_{G}\left[  r\right]  }{\lambda_{G}}=L_{P}\left[  \frac{r\left(
x\right)  }{x-\omega}\right]  +Mr\left(  \omega\right)  ,\quad r\in
\mathbb{C}\left[  x\right]  ,
\]
or by%
\[
\left(  z-\omega\right)  \frac{S_{G}\left(  z\right)  }{\lambda_{G}}=S\left(
z\right)  -S\left(  \omega\right)  +M,
\]
where $M-S\left(  \omega\right)  \neq0$ and%
\[
\lambda_{G}=\frac{L_{G}\left[  1\right]  }{M-S\left(  \omega\right)  }.
\]

\end{enumerate}

\strut If the coefficients in the 3-term recurrence relation (\ref{3-term})
are constant, then we see that $P_{n}^{\ast}\left(  x\right)  $ and
$P_{n-1}\left(  x\right)  $ satisfy the same recurrence, and have the same
initial conditions. Therefore, they are identical $P_{n}^{\ast}\left(
x\right)  =P_{n-1}\left(  x\right)  $ and we obtain%
\[
x=\frac{P_{n+1}}{P_{n}}+\beta+\gamma\frac{P_{n-1}}{P_{n}}=\frac{P_{n+1}%
}{P_{n+1}^{\ast}}+\beta+\gamma\frac{P_{n}^{\ast}}{P_{n}}.
\]
Using Markov's theorem, we conclude that%
\[
z=\frac{L_{P}\left[  1\right]  }{S_{P}\left(  z\right)  }+\beta+\gamma
\frac{S_{P}\left(  z\right)  }{L_{P}\left[  1\right]  }.
\]
Assuming that $L_{P}\left[  1\right]  =1,$ we find that $S_{P}\left(
z\right)  $ is the solution of the quadratic equation%
\begin{equation}
\gamma S_{P}^{2}-\left(  z-\beta\right)  S_{P}+1=0. \label{quad1}%
\end{equation}

Multiplying (\ref{quadratic}) by $\gamma$, (\ref{quad1}) by $AQ,$ and
subtracting, we get%
\[
\left[  \gamma\left(  BQ+uA-C\right)  +\left(  z-\beta\right)  AQ\right]
S_{P}+u\gamma B-AQ=0.
\]
Therefore,%
\[
A=u\gamma,\quad B=Q,\quad C=Q^{2}+u\left(  z-\beta\right)  Q+u^{2}\gamma,
\]
or replacing $Q\left(  z\right)  $ by $\left(  1-u\right)  z+u\beta-v,$
\[
B=\left(  1-u\right)  z+u\beta-v,\quad C=\left(  z-v\right)  \left[  \left(
1-u\right)  z+u\beta-v\right]  +u^{2}\gamma,
\]
We conclude that%
\begin{equation}
\frac{S_{q}\left(  z\right)  }{L_{q}\left[  1\right]  }=\frac{u\gamma
S_{P}\left(  z\right)  +\left(  1-u\right)  z+u\beta-v}{\left(  z-v\right)
\left[  \left(  1-u\right)  z+u\beta-v\right]  +u^{2}\gamma}.
\label{Transformed S}%
\end{equation}

Polynomial solutions of 3-term recurrence relations with constant coefficients%
\[
xq_{n}=q_{n+1}+\beta q_{n}+\gamma q_{n-1},\quad q_{0}\left(  x\right)
=u,\quad q_{1}\left(  x\right)  =x-v,
\]
were analyzed in \cite{MR766488}, where the authors concluded that the linear
functional $L_{q}$ was of the form%
\[
\frac{L_{q}\left[  r\right]  }{L_{q}\left[  f\right]  }=\frac{1}{2\pi\gamma}%
{\displaystyle\int\limits_{\beta-2\sqrt{\gamma}}^{\beta+2\sqrt{\gamma}}}
r\left(  x\right)  \frac{\sqrt{4\gamma-\left(  x-\beta\right)  ^{2}}}{f\left(
x\right)  }dx+M_{1}r\left(  y_{1}\right)  +M_{2}r\left(  y_{2}\right)  ,\quad
r\in\mathbb{C}\left[  x\right]  ,
\]
with%
\[
f\left(  x\right)  =\left(  x-v\right)  \left[  \left(  1-u\right)
x+u\beta-v\right]  +u^{2}\gamma=\left(  1-u\right)  \left(  x-y_{1}\right)
\left(  x-y_{2}\right)  .
\]
They only considered the case where $u>0$ and all parameters (including the
roots of $f)$ are real numbers. For the masses $M_{1},M_{2}$ they obtained
\[
M_{i}=\frac{2}{\sqrt{\left[  \left(  \beta-v\right)  ^{2}+4\gamma\left(
u-1\right)  \right]  _{+}}}\left[  \frac{u\gamma}{\left\vert y_{i}%
-v\right\vert }-\frac{\left\vert y_{i}-v\right\vert }{u}\right]  _{+},\quad
u\neq1,
\]
and
\[
M_{1}=M_{2}=\left[  1-\frac{\gamma}{\left(  \beta-v\right)  ^{2}}\right]
_{+},\quad u=1,
\]
where
\[
\left[  x\right]  _{+}=\frac{x+\left\vert x\right\vert }{2}.
\]

The solution of (\ref{quad1}) having the right asymptotic behavior, is
\begin{equation}
S_{P}\left(  z\right)  =\frac{z-\beta-\sqrt{\left(  z-\beta\right)
^{2}\allowbreak-4\gamma}}{2\gamma}\sim\frac{1}{z}+\frac{\beta}{z^{2}},\quad
z\rightarrow\infty, \label{Sp}%
\end{equation}
and it follows that an integral representation of the lineal functional
$L_{P}$ is given by%
\[
L_{P}\left[  r\right]  =\frac{1}{2\pi\gamma}%
{\displaystyle\int\limits_{\beta-2\sqrt{\gamma}}^{\beta+2\sqrt{\gamma}}}
r\left(  x\right)  \sqrt{4\gamma-\left(  x-\beta\right)  ^{2}}dx,\quad
r\in\mathbb{C}\left[  x\right]  .
\]
Using the change of variables
\[
x=\beta+2y\sqrt{\gamma},
\]
we obtain%
\[
L_{P}\left[  r\right]  =\frac{2}{\pi}%
{\displaystyle\int\limits_{-1}^{1}}
r\left(  \beta+2y\sqrt{\gamma}\right)  \sqrt{1-y^{2}}dy.
\]
Therefore, it is enough to study the linear functional%
\[
L_{U}\left[  r\right]  =\frac{2}{\pi}%
{\displaystyle\int\limits_{-1}^{1}}
r\left(  y\right)  \sqrt{1-y^{2}}dy,
\]
associated to the Chebyshev polynomials of the second kind, which we will
define in the next subsection.

\subsection{Chebyshev polynomials}

The monic \emph{Chebyshev polynomials of the second kind} $U_{n}\left(
x\right)  $ are defined by \cite[9.8.36]{MR2656096}%
\[
U_{n}\left(  x\right)  =2^{-n}\left(  n+1\right)  \ _{2}F_{1}\left(
\begin{array}
[c]{c}%
-n,n+2\\
\frac{3}{2}%
\end{array}
;\frac{1-x}{2}\right)  .
\]
They are a solution of the recurrence relation%
\begin{equation}
U_{n+1}-xU_{n}+\frac{1}{4}U_{n-1}=0, \label{reqcheb}%
\end{equation}
with initial conditions \cite[9.8.40]{MR2656096}%
\begin{equation}
U_{0}\left(  x\right)  =1,\quad U_{1}\left(  x\right)  =x. \label{iniU}%
\end{equation}
Note that
\begin{equation}
U_{-1}\left(  x\right)  =0. \label{Um1}%
\end{equation}

The polynomials $U_{n}\left(  x\right)  $ satisfy the orthogonality relation
\cite[9.8.38]{MR2656096}%
\begin{equation}
L_{U}\left[  U_{n}U_{m}\right]  =\frac{2}{\pi}%
{\displaystyle\int\limits_{-1}^{1}}
U_{n}\left(  x\right)  U_{m}\left(  x\right)  \sqrt{1-x^{2}}dx=\delta_{n,m}.
\label{orthoU}%
\end{equation}
The Stieltjes transforms of the linear functional $L_{U}$ is given by
\cite{MR903848}%
\begin{equation}
S_{U}\left(  z\right)  =\frac{2}{\pi}%
{\displaystyle\int\limits_{-1}^{1}}
\frac{\sqrt{1-x^{2}}}{z-x}dx=2z\left(  1-\sqrt{1-z^{-2}}\right)  ,\quad
z\in\mathbb{C}\setminus\left[  -1,1\right]  , \label{SU}%
\end{equation}
where here and in the rest of the paper
\[
\sqrt{\quad}:\mathbb{C}\rightarrow\left\{  z\in\mathbb{C}\mid-\frac{\pi}%
{2}<\arg\left(  z\right)  \leq\frac{\pi}{2}\right\}
\]
denotes the principal branch of the square root. Note that%
\[
\sqrt{1-z^{-2}}\sim1,\quad z\rightarrow\infty.
\]

It is clear from the three-term recurrence relation (\ref{Reqn}) that the
polynomials $D_{n,k}(x;a)$ are related to the polynomials $U_{n}\left(
x\right)  $. Let's introduce the scaled polynomials $d_{n}(x;k)$ defined by%
\begin{equation}
d_{n}(x;k)=\left(  2\sqrt{a}\right)  ^{-n}D_{n,k}\left(  2\sqrt{a}x\right)  .
\label{d-scaled}%
\end{equation}
The polynomials $d_{n}(x;k)$ are a solution of the recurrence (\ref{reqcheb})
satisfied by the polynomials $U_{n}\left(  x\right)  $, with initial
conditions
\[
d_{0}(x;k)=2-k,\quad d_{1}(x;k)=x.
\]
Therefore, we see that the scaled polynomials $d_{n}(x;k)$ are co-recursive
polynomials (with respect to the Chebyshev polynomials of the second kind)
with parameters $u=2-k,$ $v=0.$

Using (\ref{q1}), we have%
\begin{equation}
d_{n}(x;k)=\left(  2-k\right)  U_{n}\left(  x\right)  +\left(  k-1\right)
xU_{n-1}\left(  x\right)  , \label{Cheby}%
\end{equation}
since%
\[
d_{n}^{\ast}(x;k)=U_{n}^{\ast}\left(  x\right)  =U_{n-1}\left(  x\right)  .
\]

\begin{remark}
If we use the values of the monic Chebyshev polynomials at $x=0$
\cite[18.6.1]{MR2723248}%
\[
U_{n}\left(  0\right)  =2^{-n}\cos\left(  \frac{n\pi}{2}\right)  ,
\]
and the representation (\ref{Cheby}), we get%
\[
D_{n,k}(0;a)=\left(  2\sqrt{a}\right)  ^{n}d_{n}(0;k)=\left(  \sqrt{a}\right)
^{n}\left(  2-k\right)  \cos\left(  \frac{n\pi}{2}\right)  ,
\]
in agreement with (\ref{even})-(\ref{odd}).
\end{remark}

If $k=0,1,2,$ (\ref{d-scaled}) and (\ref{Cheby}) give%
\[
D_{n,k}\left(  x\right)  =\left(  2\sqrt{a}\right)  ^{n}\left[  \left(
2-k\right)  U_{n}\left(  \frac{x}{2\sqrt{a}}\right)  +\left(  k-1\right)
\frac{x}{2\sqrt{a}}U_{n-1}\left(  \frac{x}{2\sqrt{a}}\right)  \right]
\]%
\begin{align*}
D_{n,0}(x;a)  &  =\left(  2\sqrt{a}\right)  ^{n}\left[  2U_{n}\left(  \frac
{x}{2\sqrt{a}}\right)  -\frac{x}{2\sqrt{a}}U_{n-1}\left(  \frac{x}{2\sqrt{a}%
}\right)  \right]  ,\\
D_{n,1}(x;a)  &  =\left(  2\sqrt{a}\right)  ^{n}U_{n}\left(  \frac{x}%
{2\sqrt{a}}\right)  ,
\end{align*}%
\begin{equation}
D_{n,2}(x;a)=\left(  2\sqrt{a}\right)  ^{n-1}xU_{n-1}\left(  \frac{x}%
{2\sqrt{a}}\right)  , \label{D2}%
\end{equation}
and in particular, for $a=1,$ we have%
\begin{align*}
D_{n}(2x)  &  =D_{n,0}(2x;1)=2^{n+1}U_{n}\left(  x\right)  -x2^{n}%
U_{n-1}\left(  x\right)  ,\\
E_{n}(2x)  &  =D_{n,1}(2x;1)=2^{n}U_{n}\left(  x\right)  ,
\end{align*}
as it was observed in \cite{MR2928474}.

\begin{remark}
Stoll \cite{MR2398338} studied second order recurrences with constant
coefficients and general initial conditions and found polynomial
decompositions in terms of Chebyshev polynomials.
\end{remark}

\section{Main results \label{Section3}}

In this section we find a representation for the linear functional $L_{k}$
defined by (\ref{L1}). Although this seems to be something already considered
in \cite{MR766488}, we have found that the range of parameters of the
polynomials $D_{n,k}\left(  x\right)  $ requires a different analysis.

Let $k\neq2$, $l_{k}\left[  r\right]  $ denote the linear functional
satisfying
\begin{equation}
l_{k}\left[  1\right]  =\frac{1}{2-k},\quad l_{k}\left[  d_{n}d_{m}\right]
=4^{-n}\delta_{n,m},\quad n,m\in\mathbb{N}, \label{l}%
\end{equation}
and $s(z;k)$ be it's Stieltjes transform. From (\ref{Transformed S}) and
(\ref{SU}), we have
\begin{equation}
s(z;k)=\frac{2z}{2-k}\frac{\left(  k-2\right)  \sqrt{1-z^{-2}}+k}{4\left(
k-1\right)  z^{2}+\left(  k-2\right)  ^{2}}, \label{S}%
\end{equation}
since $v=\beta=0.$

Our next objective is to represent the function $s(z;k)$ as the Stieltjes
transform of a distribution. We begin with a couple of lemmas.

\begin{lemma}
\label{Lemma1}Let $\mathcal{L}$ be a linear functional with Stieltjes
transform $S\left(  z\right)  $ and%
\[
f\left(  x\right)  =\left(  x-\omega_{1}\right)  \left(  x-\omega_{2}\right)
.
\]
Then,%
\[
\mathcal{L}\left[  \frac{1}{f\left(  x\right)  \left(  z-x\right)  }\right]
=\frac{S\left(  z\right)  }{f\left(  z\right)  }+\frac{\left[  S\left(
\omega_{2}\right)  -S\left(  \omega_{1}\right)  \right]  z+\omega_{2}S\left(
\omega_{1}\right)  -\omega_{1}S\left(  \omega_{2}\right)  }{\left(  \omega
_{1}-\omega_{2}\right)  f\left(  z\right)  },
\]
where we always assume that the functional $\mathcal{L}$ acts on the variable
$x.$
\end{lemma}

\begin{proof}
Since%
\begin{gather*}
\frac{1}{\left(  x-\omega_{1}\right)  \left(  x-\omega_{2}\right)  \left(
z-x\right)  }\allowbreak=\frac{1}{\left(  z-\omega_{1}\right)  \left(
z-\omega_{2}\right)  \left(  z-x\right)  }\\
+\frac{1}{\left(  \omega_{1}-\omega_{2}\right)  \left(  z-\omega_{2}\right)
\left(  \omega_{2}-x\right)  }-\frac{1}{\left(  \omega_{1}-\omega_{2}\right)
\left(  z-\omega_{1}\right)  \left(  \omega_{1}-x\right)  },
\end{gather*}
we get%
\[
\mathcal{L}\left[  \frac{1}{f\left(  x\right)  \left(  z-x\right)  }\right]
=\frac{S\left(  z\right)  }{\left(  z-\omega_{1}\right)  \left(  z-\omega
_{2}\right)  }+\frac{S\left(  \omega_{2}\right)  }{\left(  \omega_{1}%
-\omega_{2}\right)  \left(  z-\omega_{2}\right)  }-\frac{S\left(  \omega
_{1}\right)  }{\left(  \omega_{1}-\omega_{2}\right)  \left(  z-\omega
_{1}\right)  },
\]
and the conclusion follows.
\end{proof}

For the particular case of the linear functional $L_{U}$ defined by
(\ref{orthoU}) with Stieltjes transform $S_{U}\left(  z\right)  $ given in
(\ref{SU}) and $f\left(  x\right)  =x^{2}-b^{2}$, Lemma \ref{Lemma1} gives%
\begin{equation}
L_{U}\left[  \frac{1}{\left(  x^{2}-b^{2}\right)  \left(  z-x\right)
}\right]  =2z\frac{\sqrt{1-b^{-2}}-\sqrt{1-z^{-2}}}{z^{2}-b^{2}},\quad
z,b\in\mathbb{C}\setminus\left[  -1,1\right]  , \label{Integral}%
\end{equation}
because%
\begin{gather*}
\frac{S_{U}\left(  z\right)  }{z^{2}-b^{2}}-\frac{S_{U}\left(  b\right)
}{2b\left(  z-b\right)  }+\frac{S_{U}\left(  -b\right)  }{2b\left(
z+b\right)  }=\frac{2z\left(  1-\sqrt{1-z^{-2}}\right)  }{z^{2}-b^{2}}\\
+\frac{2b\left(  1-\sqrt{1-b^{-2}}\right)  }{2b\left(  b-z\right)  }%
+\frac{2\left(  -b\right)  \left(  1-\sqrt{1-b^{-2}}\right)  }{2b\left(
b+z\right)  }.
\end{gather*}

Note that since \cite{MR903848}%
\begin{equation}
\frac{2}{\pi}%
{\displaystyle\int\limits_{-1}^{1}}
\frac{1}{\sqrt{1-x^{2}}}\frac{dx}{z-x}=\frac{2}{z\sqrt{1-z^{-2}}},\quad
z\in\mathbb{C}\setminus\left[  -1,1\right]  , \label{ST}%
\end{equation}
we have%
\begin{gather*}
\underset{b^{2}\rightarrow1^{-}}{\lim}L_{U}\left[  \frac{1}{\left(
x^{2}-b^{2}\right)  \left(  z-x\right)  }\right]  =\frac{2}{\pi}%
{\displaystyle\int\limits_{-1}^{1}}
\frac{\sqrt{1-x^{2}}}{x^{2}-1}\frac{1}{z-x}dx\\
=-\frac{2}{\pi}%
{\displaystyle\int\limits_{-1}^{1}}
\frac{1}{\sqrt{1-x^{2}}\left(  z-x\right)  }dx=\frac{-2}{z\sqrt{1-z^{-2}}}\\
=-2z\frac{\sqrt{1-z^{-2}}}{z^{2}-1}=\underset{b^{2}\rightarrow1^{-}}{\lim
}2z\frac{\sqrt{1-b^{-2}}-\sqrt{1-z^{-2}}}{z^{2}-b^{2}},
\end{gather*}
and in particular%
\begin{equation}
\frac{1}{2z\sqrt{1-z^{-2}}}=\frac{1}{4}L_{U}\left[  \frac{1}{\left(
1-x^{2}\right)  \left(  z-x\right)  }\right]  . \label{b=1}%
\end{equation}

\begin{lemma}
\label{Lemma2}Let $\omega\left(  k\right)  $ be defined by
\begin{equation}
\omega\left(  k\right)  =\frac{1}{2}\frac{k-2}{\sqrt{k-1}}\mathrm{i},\quad
k\neq1, \label{w}%
\end{equation}
where $\mathrm{i}^{2}=-1.$ Then,
\[
\omega\left(  k\right)  \in\mathbb{C}\setminus\left[  -1,1\right]  ,\quad
k\in\mathbb{R}\setminus\left\{  0,1,2\right\}  .
\]

\end{lemma}

\begin{proof}
The result follows immediately from the definition (\ref{w}), since%
\begin{align*}
\omega\left(  k\right)   &  \in\left(  -\infty,-1\right)  ,\quad k\in\left(
-\infty,0\right)  \cup\left(  0,1\right)  ,\\
\mathrm{i}\omega\left(  k\right)   &  \in\mathbb{R},\quad k\in\left(
1,\infty\right)  ,\\
\omega\left(  0\right)   &  =-1,\quad\omega\left(  2\right)  =0.
\end{align*}

\end{proof}

The function $s(z;k)$ defined in (\ref{S}) has a branch cut on the segment
$\left[  -1,1\right]  $ and (perhaps removable) poles at $z=\pm\omega$ if
$k\neq1.$ In the next theorem we split $s(z;k)$ in two parts, one analytic in
$\mathbb{C}\setminus\left[  -1,1\right]  $ and the other analytic in
$\mathbb{C}\setminus\left\{  \pm\omega\right\}  .$

\begin{theorem}
\label{Theorem}Let $k\neq2$ and $z\in\mathbb{C}\setminus\left[  -1,1\right]
.$ Then, we have%
\begin{equation}
s(z;k)=s_{c}(z;k)+\chi\left(  k\right)  s_{d}(z;k), \label{s1}%
\end{equation}
where $\chi\left(  k\right)  $ is the characteristic function defined by%
\begin{equation}
\chi\left(  k\right)  =\left\{
\begin{array}
[c]{c}%
0,\quad k\in\left[  0,2\right] \\
1,\quad k\in\mathbb{R}\setminus\left[  0,2\right]
\end{array}
\right.  , \label{char}%
\end{equation}
$s_{c}(z;k)$ is the continuous part of $s(z;k)$
\begin{equation}
s_{c}(z;k)=\frac{2}{\pi}%
{\displaystyle\int\limits_{-1}^{1}}
\frac{\sqrt{1-x^{2}}}{4\left(  k-1\right)  x^{2}+\left(  k-2\right)  ^{2}%
}\frac{1}{z-x}dx, \label{sc}%
\end{equation}
and $s_{d}(z;k)$ is the discrete part of $s(z;k)$
\begin{equation}
s_{d}\left(  z;k\right)  =\frac{4k}{2-k}\frac{z}{4\left(  k-1\right)
z^{2}+\left(  k-2\right)  ^{2}}. \label{sd}%
\end{equation}

\end{theorem}

\begin{proof}
Let $k\in\mathbb{R}\setminus\left\{  0,1,2\right\}  .$ From (\ref{S}), we have%
\[
s(z;k)=2z\frac{\frac{k}{2-k}-\sqrt{1-z^{-2}}}{4\left(  k-1\right)
z^{2}+\left(  k-2\right)  ^{2}}.
\]
Using (\ref{w}), we get%
\begin{equation}
\frac{\left(  k-2\right)  ^{2}}{4\left(  k-1\right)  }=-\omega^{2}%
,\quad1-\omega^{-2}\left(  k\right)  =\left(  \frac{k}{2-k}\right)  ^{2}.
\label{1-w}%
\end{equation}
Hence,
\[
s(z;k)=\frac{2z}{4\left(  k-1\right)  }\frac{\sqrt{1-\omega^{-2}}%
-\sqrt{1-z^{-2}}+\frac{k}{2-k}-\left\vert \frac{k}{2-k}\right\vert }%
{z^{2}-\omega^{2}}.
\]

Since we know from Lemma \ref{Lemma2} that $\omega\left(  k\right)
\in\mathbb{C}\setminus\left[  -1,1\right]  $, we can use (\ref{Integral}) with
$b=\omega$ and obtain%
\[
s(z;k)=\frac{1}{4\left(  k-1\right)  }L_{U}\left[  \frac{1}{\left(
x^{2}-\omega^{2}\right)  \left(  z-x\right)  }\right]  +\frac{2z}{4\left(
k-1\right)  }\frac{\frac{k}{2-k}-\left\vert \frac{k}{2-k}\right\vert }%
{z^{2}-\omega^{2}}.
\]
But%
\[
\frac{k}{2-k}-\left\vert \frac{k}{2-k}\right\vert =\left\{
\begin{array}
[c]{c}%
0,\quad k\in\lbrack0,2)\\
\frac{2k}{2-k},\quad k\in\mathbb{R}\setminus\left[  0,2\right]
\end{array}
\right.  ,
\]
and therefore%
\[
s(z;k)=\left\{
\begin{array}
[c]{c}%
s_{c}(z;k),\quad k\in\left(  0,2\right)  \setminus\left\{  1\right\} \\
s_{c}(z;k)+s_{d}(z;k),\quad k\in\mathbb{R}\setminus\left[  0,2\right]
\end{array}
\right.  ,
\]
where%
\[
s_{c}(z;k)=L_{U}\left[  \frac{1}{4\left(  k-1\right)  x^{2}+\left(
k-2\right)  ^{2}}\frac{1}{z-x}\right]  ,
\]
and
\[
s_{d}(z;k)=\frac{2z}{4\left(  k-1\right)  }\frac{2k}{2-k}\frac{1}{z^{2}%
-\omega^{2}}.
\]

If $k=0,$ then we have from (\ref{S})%
\[
s(z;0)=z\frac{\sqrt{1-z^{-2}}}{2\left(  z^{2}-1\right)  }=\frac{1}%
{2z\sqrt{1-z^{-2}}},
\]
and using (\ref{b=1})$,$ we get%
\[
\frac{1}{2z\sqrt{1-z^{-2}}}=\frac{1}{4}L_{U}\left[  \frac{1}{\left(
1-x^{2}\right)  \left(  z-x\right)  }\right]  =s_{c}(z;0).
\]

If $k=1,$ then we have from (\ref{S})%
\[
s(z;1)=2z\left(  1-\sqrt{1-z^{-2}}\right)  ,
\]
and using (\ref{SU}) we get%
\[
s(z;1)=\frac{2}{\pi}%
{\displaystyle\int\limits_{-1}^{1}}
\frac{\sqrt{1-x^{2}}}{z-x}dx=s_{c}(z;1).
\]

\end{proof}

\begin{remark}
Orthogonal polynomials with linear functionals of the form
\[
\mathcal{L}\left[  r\right]  =%
{\displaystyle\int\limits_{-1}^{1}}
r\left(  x\right)  \frac{\left(  1-x\right)  ^{\alpha}\left(  1+x\right)
^{\beta}}{f\left(  x\right)  }dx,\quad\alpha,\beta=\pm\frac{1}{2},
\]
where $f(x)$ is a polynomial, are called \emph{Bernstein--Szeg\"{o}
polynomials} \cite[18.31]{MR2723248}, \cite[2.6]{MR0372517}. These polynomials
are examples of the Geronimus transformation applied to the Jacobi
polynomials, see \cite{MR2770493}, \cite{MR1502972}, \cite[2.7.3]{MR2191786},
and \cite{MR1482157}.
\end{remark}

\begin{corollary}
Let $k\neq2,$ the linear functional $l_{k}$ be defined by (\ref{l}), the
characteristic function $\chi\left(  k\right)  $ be defined by (\ref{char}),
and the function $\omega\left(  k\right)  $ de defined by \ref{w}.

Then, for all $r\left(  x\right)  \in\mathbb{C}\left[  x\right]  $ we have%
\begin{equation}
l_{k}\left[  r\right]  =l_{k}^{\left(  c\right)  }\left[  r\right]
+\chi\left(  k\right)  l_{k}^{\left(  d\right)  }\left[  r\right]  ,
\label{l rep}%
\end{equation}
where%
\begin{equation}
l_{k}^{\left(  c\right)  }\left[  r\right]  =\frac{2}{\pi}%
{\displaystyle\int\limits_{-1}^{1}}
\frac{r\left(  x\right)  \sqrt{1-x^{2}}}{4\left(  k-1\right)  x^{2}+\left(
k-2\right)  ^{2}}dx, \label{lc}%
\end{equation}%
\begin{equation}
l_{k}^{\left(  d\right)  }\left[  r\right]  =k\frac{r\left(  \omega\right)
+r\left(  -\omega\right)  }{2\left(  k-1\right)  \left(  2-k\right)  },
\label{ld}%
\end{equation}
and we assume that%
\begin{equation}
\underset{k\rightarrow1}{\lim}\frac{\chi\left(  k\right)  }{k-1}=0.
\label{limits}%
\end{equation}

\end{corollary}

\begin{proof}
The result is a direct consequence of Theorem \ref{Theorem}, since%
\begin{gather*}
\frac{4k}{2-k}\frac{z}{4\left(  k-1\right)  z^{2}+\left(  k-2\right)  ^{2}%
}=\frac{k}{\left(  2-k\right)  \left(  k-1\right)  }\frac{z}{z^{2}-\omega^{2}%
}\\
=\frac{1}{2}\frac{k}{\left(  k-1\right)  \left(  2-k\right)  }\left(  \frac
{1}{z-\omega}+\frac{1}{z+\omega}\right)  .
\end{gather*}

For $k=1,$ we see from (\ref{Cheby}) that%
\[
d_{n}(x;1)=U_{n}\left(  x\right)  ,
\]
and therefore%
\[
l_{1}\left[  r\right]  =L_{U}\left[  r\right]  =\frac{2}{\pi}%
{\displaystyle\int\limits_{-1}^{1}}
r\left(  x\right)  \sqrt{1-x^{2}}dx,
\]
which agrees with (\ref{l rep}) for $k=1$ if we use (\ref{limits}).
\end{proof}

\begin{remark}
For $k=2,$ we see from (\ref{Cheby}) that%
\[
d_{n}(x;2)=xU_{n-1}\left(  x\right)  ,
\]
and therefore we can interpret $l_{2}$ as the linear functional%
\[
l_{2}\left[  r\right]  =L_{U}\left[  \frac{r\left(  x\right)  }{4x^{2}%
}\right]  =\frac{1}{2\pi}%
{\displaystyle\int\limits_{-1}^{1}}
r\left(  x\right)  \frac{\sqrt{1-x^{2}}}{x^{2}}dx,
\]
defined for all polynomials $r\left(  x\right)  $ such that $r\left(
x\right)  =x^{2}p\left(  x\right)  ,$ $p\left(  x\right)  \in\mathbb{C}\left[
x\right]  .$ It follows that $d_{n}(x;2)$ will be a family of orthogonal
polynomials for $n\geq1.$
\end{remark}

\subsection{The Dickson polynomials}

We can now use the previous results to the Dickson polynomials of the
$(k+1)$-th kind $D_{n,k}(x;a),$ related to the scaled polynomials $d_{n}(x;k)$
by (\ref{d-scaled}).

\begin{lemma}
Let $\mathfrak{L}_{k}:\mathbb{C}\left[  x\right]  \rightarrow\mathbb{C}$ be
the linear functional defined by%
\begin{equation}
\mathfrak{L}_{k}\left[  r\left(  x\right)  \right]  =l_{k}\left[  r\left(
2\sqrt{a}x\right)  \right]  , \label{L-rep}%
\end{equation}
where $l_{k}$ is the linear functional defined by (\ref{l}). Then,
$\mathfrak{L}_{k}$ satisfies
\begin{equation}
\mathfrak{L}_{k}\left[  D_{0,k}^{2}\right]  =2-k,\quad\mathfrak{L}_{k}\left[
D_{n,k}D_{m,k}\right]  =a^{n}\delta_{n,m},\quad n,m\in\mathbb{N},
\label{orthoLk}%
\end{equation}
and for $k\neq2,$ its Stieltjes transform is given by%
\begin{equation}
S\left(  z;k,a\right)  =\frac{z}{2\left(  2-k\right)  }\frac{\left(
k-2\right)  \sqrt{1-4az^{-2}}+k}{\left(  k-1\right)  z^{2}+\left(  k-2\right)
^{2}a},\quad z\in\mathbb{C}\setminus\left[  -2\sqrt{a},2\sqrt{a}\right]  .
\label{SD}%
\end{equation}

\end{lemma}

\begin{proof}
Using (\ref{d-scaled}) in (\ref{L-rep}), we have%
\begin{align*}
\mathfrak{L}_{k}\left[  D_{n,k}\left(  x\right)  D_{m,k}\left(  x\right)
\right]   &  =l_{k}\left[  D_{n,k}\left(  2\sqrt{a}x\right)  D_{m,k}\left(
2\sqrt{a}x\right)  \right] \\
&  =\left(  2\sqrt{a}\right)  ^{n+m}l_{k}\left[  d_{n}(x;k)d_{m}(x;k)\right]
.
\end{align*}
Therefore, (\ref{l}) gives
\[
\mathfrak{L}_{k}\left[  D_{0,k}^{2}\right]  =l_{k}\left[  \left(  2-k\right)
^{2}\right]  =\left(  2-k\right)  ^{2}l_{k}\left[  1\right]  =2-k,
\]
and%
\[
\mathfrak{L}_{k}\left[  D_{n,k}\left(  x\right)  D_{m,k}\left(  x\right)
\right]  =\left(  2\sqrt{a}\right)  ^{n+m}4^{-n}\delta_{n,m},\quad
n,m\in\mathbb{N}.
\]
But since%
\[
\left(  2\sqrt{a}\right)  ^{n+m}4^{-n}\delta_{n,m}=a^{n}\delta_{n,m},
\]
(\ref{orthoLk}) follows.

Using (\ref{S}) in (\ref{L-rep}), we get%
\begin{gather*}
\mathfrak{L}_{k}\left[  \frac{1}{z-x}\right]  =l_{k}\left[  \frac{1}%
{z-2\sqrt{a}x}\right]  =\frac{1}{2\sqrt{a}}l_{k}\left[  \frac{1}{\frac
{z}{2\sqrt{a}}-x}\right] \\
=\frac{1}{2\sqrt{a}}s\left(  \frac{z}{2\sqrt{a}};k\right)  =\frac{1}{2\sqrt
{a}}\frac{2\frac{z}{2\sqrt{a}}}{2-k}\frac{\left(  k-2\right)  \sqrt
{1-4az^{-2}}+k}{4\left(  k-1\right)  \frac{z^{2}}{4a}+\left(  k-2\right)
^{2}},
\end{gather*}
$\allowbreak$and (\ref{SD}) follows.
\end{proof}

\begin{corollary}
The linear functional $\mathfrak{L}_{k}$ defined by (\ref{L-rep}) is identical
to the linear functional $L_{k}$ satisfying (\ref{L1}).
\end{corollary}

Next, we find a representation for the linear functional $L_{k}.$

\begin{theorem}
\label{Main}Let $k\neq2$ and $L_{k}$ be the linear functional defined by
(\ref{L1}). Then, $L_{k}$ admits the representation%
\[
L_{k}\left[  r\right]  =L_{k}^{\left(  c\right)  }\left[  r\right]
+\chi\left(  k\right)  L_{k}^{\left(  d\right)  }\left[  r\right]  ,
\]
where $\chi\left(  k\right)  $ was defined in (\ref{char}),%
\begin{equation}
L_{k}^{\left(  c\right)  }\left[  r\right]  =\frac{1}{2\pi}%
{\displaystyle\int\limits_{-2\sqrt{a}}^{2\sqrt{a}}}
\frac{r\left(  t\right)  \sqrt{4a-t^{2}}}{\left(  k-1\right)  t^{2}+\left(
k-2\right)  ^{2}a}dt \label{Lc}%
\end{equation}
and%
\begin{equation}
L_{k}^{\left(  d\right)  }\left[  r\right]  =k\frac{r\left(  \Omega\right)
+r\left(  -\Omega\right)  }{2\left(  k-1\right)  \left(  2-k\right)  },
\label{Ld}%
\end{equation}
with $\Omega\left(  k\right)  $ defined by
\begin{equation}
\Omega\left(  k\right)  =2\sqrt{a}\omega\left(  k\right)  =\sqrt{a}\frac
{k-2}{\sqrt{k-1}}\mathrm{i},\quad k\neq1. \label{Omega}%
\end{equation}

\end{theorem}

\begin{proof}
Changing variables to $t=2\sqrt{a}x$ in the integral
\[
I=%
{\displaystyle\int\limits_{-1}^{1}}
\frac{r\left(  x\right)  \sqrt{1-x^{2}}}{4\left(  k-1\right)  x^{2}+\left(
k-2\right)  ^{2}}dx,
\]
we obtain%
\[
I=\frac{1}{4}%
{\displaystyle\int\limits_{-2\sqrt{a}}^{2\sqrt{a}}}
\frac{r\left(  \frac{t}{2\sqrt{a}}\right)  \sqrt{4a-t^{2}}}{\left(
k-1\right)  t^{2}+\left(  k-2\right)  ^{2}a}dt.
\]
Therefore, from (\ref{lc}) we get%
\[
l_{k}^{\left(  c\right)  }\left[  r\left(  2\sqrt{a}t\right)  \right]
=\frac{1}{2\pi}%
{\displaystyle\int\limits_{-2\sqrt{a}}^{2\sqrt{a}}}
\frac{r\left(  t\right)  \sqrt{4a-t^{2}}}{\left(  k-1\right)  t^{2}+\left(
k-2\right)  ^{2}a}dt.
\]
Also, from (\ref{ld}) and (\ref{Omega}) we have%
\[
l_{k}^{\left(  d\right)  }\left[  r\left(  2\sqrt{a}t\right)  \right]
=k\frac{r\left(  2\sqrt{a}\omega\right)  +r\left(  -2\sqrt{a}\omega\right)
}{2\left(  k-1\right)  \left(  2-k\right)  }=k\frac{r\left(  \Omega\right)
+r\left(  -\Omega\right)  }{2\left(  k-1\right)  \left(  2-k\right)  }.
\]

Therefore, using (\ref{l rep}) in (\ref{L-rep}) we see that%
\[
L_{k}\left[  r\right]  =\frac{1}{2\pi}%
{\displaystyle\int\limits_{-2\sqrt{a}}^{2\sqrt{a}}}
\frac{r\left(  t\right)  \sqrt{4a-t^{2}}}{\left(  k-1\right)  t^{2}+\left(
k-2\right)  ^{2}a}dt+\chi\left(  k\right)  k\frac{r\left(  \Omega\right)
+r\left(  -\Omega\right)  }{2\left(  k-1\right)  \left(  2-k\right)  }.
\]

\end{proof}

Although Theorem \ref{Main} seems to be valid only when $k\neq2,$ we can see
that $L_{2}$ is well defined$.$

\begin{lemma}
Let $k\neq1$ and $\Omega\left(  k\right)  $ be defined by (\ref{Omega}). Then,%
\begin{equation}
D_{n,k}\left(  \Omega;a\right)  =\left(  2-k\right)  \left(  -\mathrm{i}%
\sqrt{\frac{a}{k-1}}\right)  ^{n}. \label{DW}%
\end{equation}

\end{lemma}

\begin{proof}
Lets assume that%
\[
D_{n,k}\left(  \Omega;a\right)  =b_{0}B^{n},
\]
for some functions $b_{0}\left(  k,a\right)  $ and $B\left(  k,a\right)  .$
Using (\ref{Reqn}), we have%
\[
0=b_{0}B^{n+2}-\Omega b_{0}B^{n+1}+ab_{0}B^{n}=b_{0}B^{n}\left(  B^{2}-\Omega
B+a\right)  .
\]

Using (\ref{Initial}), we get
\[
b_{0}=D_{0,k}(\omega;a)=2-k
\]
and%
\[
\Omega=D_{1,k}\left(  \Omega;a\right)  =\left(  2-k\right)  B\left(
k,a\right)  .
\]
Thus,%
\[
B\left(  k,a\right)  =\frac{\Omega}{2-k}=-\sqrt{\frac{a}{k-1}}\mathrm{i,}%
\]
and clearly
\[
B^{2}-\Omega B+a=0.
\]

\end{proof}

It follows from the previous Lemma that $L_{2}^{\left(  d\right)  }$ is well
defined$.$

\begin{proposition}
Let $\Omega\left(  k\right)  $ be defined by (\ref{Omega}) and $L_{k}^{\left(
d\right)  }$ be defined by (\ref{Ld}). Then, for $k\neq1$%
\[
L_{k}^{\left(  d\right)  }\left[  D_{n,k}D_{m,k}\right]  =\frac{\left(
2-k\right)  k}{k-1}\left[  \frac{1+\left(  -1\right)  ^{n+m}}{2}\right]
\left(  \mathrm{i}\sqrt{\frac{a}{k-1}}\right)  ^{n+m}.
\]

\end{proposition}

\begin{proof}
From (\ref{DW}) we have%
\[
D_{n,k}\left(  \Omega;a\right)  D_{m,k}\left(  \Omega;a\right)  =\left(
2-k\right)  ^{2}\left(  -\mathrm{i}\sqrt{\frac{a}{k-1}}\right)  ^{n+m}.
\]
Using (\ref{oddeven}), we get%
\begin{align*}
&  D_{n,k}\left(  \Omega;a\right)  D_{m,k}\left(  \Omega;a\right)
+D_{n,k}\left(  -\Omega;a\right)  D_{m,k}\left(  -\Omega;a\right) \\
&  =\left(  2-k\right)  ^{2}\left(  \mathrm{i}\sqrt{\frac{a}{k-1}}\right)
^{n+m}\left[  1+\left(  -1\right)  ^{n+m}\right]  .
\end{align*}
Thus,%
\[
L_{k}^{\left(  d\right)  }\left[  D_{n,k}D_{m,k}\right]  =\frac{\left(
2-k\right)  k}{k-1}\left[  \frac{1+\left(  -1\right)  ^{n+m}}{2}\right]
\left(  \mathrm{i}\sqrt{\frac{a}{k-1}}\right)  ^{n+m}.
\]

\end{proof}

We can now extend Theorem \ref{Main} to all values of $k.$

\begin{corollary}
Let $h_{n}\left(  k\right)  $ be defined by (\ref{hD}) and $\chi\left(
k\right)  $ be defined by (\ref{char}). Then,%
\begin{gather}
\frac{1}{2\pi}%
{\displaystyle\int\limits_{-2\sqrt{a}}^{2\sqrt{a}}}
\frac{\sqrt{4a-t^{2}}D_{n,k}\left(  t\right)  D_{m,k}\left(  t\right)
}{\left(  k-1\right)  t^{2}+a\left(  k-2\right)  ^{2}}dt\label{orthoD}\\
+\chi\left(  k\right)  \frac{\left(  2-k\right)  k}{k-1}\left[  \frac
{1+\left(  -1\right)  ^{n+m}}{2}\right]  \left(  \mathrm{i}\sqrt{\frac{a}%
{k-1}}\right)  ^{n+m}=h_{n}\left(  k\right)  \delta_{n,m},\quad n,m\in
\mathbb{N}_{0}.\nonumber
\end{gather}

\end{corollary}

\begin{remark}
\label{Remark3}If we set $k=2$ in (\ref{orthoD}), we obtain%
\begin{equation}
\frac{1}{2\pi}%
{\displaystyle\int\limits_{-2\sqrt{a}}^{2\sqrt{a}}}
\frac{\sqrt{4a-t^{2}}}{t^{2}}D_{n,2}\left(  t\right)  D_{m,2}\left(  t\right)
dt=h_{n}\left(  k\right)  \delta_{n,m}, \label{inte1}%
\end{equation}
which seems to make no sense, since the integrand is singular at $t=0.$
However, if we use (\ref{D2}) we have%
\[
D_{n,2}(t;a)=a^{\frac{1}{2}\left(  n-1\right)  }tU_{n-1}\left(  \frac
{t}{2\sqrt{a}}\right)  ,
\]
and we can write (\ref{inte1}) as%
\[
a^{\frac{1}{2}\left(  n+m\right)  -1}\frac{1}{2\pi}%
{\displaystyle\int\limits_{-2\sqrt{a}}^{2\sqrt{a}}}
\sqrt{4a-t^{2}}U_{n-1}\left(  \frac{t}{2\sqrt{a}}\right)  U_{m-1}\left(
\frac{t}{2\sqrt{a}}\right)  dt=h_{n}\left(  k\right)  \delta_{n,m},
\]
or changing variables to $\tau=\frac{t}{2\sqrt{a}}$%
\[
a^{\frac{n+m}{2}}\frac{2}{\pi}%
{\displaystyle\int\limits_{-1}^{1}}
\sqrt{1-\tau^{2}}U_{n-1}\left(  \tau\right)  U_{m-1}\left(  \tau\right)
d\tau=h_{n}\left(  k\right)  \delta_{n,m}.
\]
This agrees with (\ref{orthoU}), since we have $U_{-1}=0$ from (\ref{Um1}) and
$h_{0}\left(  2\right)  =0$ from (\ref{hD}), while for $n,m\geq1$ we know from
(\ref{hD}) that%
\[
a^{-\frac{n+m}{2}}h_{n}\left(  k\right)  \delta_{n,m}=\delta_{n,m}.
\]

\end{remark}

Finally, we will use the function $S\left(  z;k,a\right)  $ to find explicit
expressions for the moments of $L_{k}.$

\begin{proposition}
Let $L_{k}$ be the linear functional defined by (\ref{L1}). Then, the moments
of $L_{k}$ of even order
\[
\mu_{2n}\left(  k\right)  =L_{k}\left[  x^{2n}\right]  ,
\]
are given by%
\begin{equation}
\mu_{2n}\left(  1\right)  =2^{2n+1}\binom{\frac{1}{2}}{n+1}\left(  -a\right)
^{n},\quad n=0,1,\ldots, \label{mu1}%
\end{equation}%
\begin{equation}
\mu_{2n}\left(  2\right)  =-2^{2n-1}\binom{\frac{1}{2}}{n}\left(  -a\right)
^{n},\quad n=1,2,\ldots, \label{mu2}%
\end{equation}
and if $k\neq1,2,$
\begin{equation}
\mu_{2n}\left(  k\right)  =-\frac{1}{2}\frac{\left(  k-2\right)  ^{2n}%
}{\left(  k-1\right)  ^{n+1}}\left(  -a\right)  ^{n}\left(  \frac{k}{k-2}+%
{\displaystyle\sum\limits_{j=0}^{n}}
\binom{\frac{1}{2}}{j}\left[  \frac{4\left(  k-1\right)  }{\left(  k-2\right)
^{2}}\right]  ^{j}\right)  . \label{moments}%
\end{equation}

\end{proposition}

\begin{proof}
From (\ref{asympt}) and (\ref{SD}), we have%
\[%
{\displaystyle\sum\limits_{j=0}^{\infty}}
\frac{\mu_{j}\left(  k\right)  }{z^{j+1}}=-\frac{z}{2}\frac{\sqrt{1-4az^{-2}%
}+\frac{k}{k-2}}{\left(  k-1\right)  z^{2}+a\left(  k-2\right)  ^{2}}.
\]
Using (\ref{oddmu}), we get%
\[%
{\displaystyle\sum\limits_{j=0}^{\infty}}
\frac{\mu_{2j}\left(  k\right)  }{z^{2j}}=-\frac{1}{2}\frac{\sqrt{1-4az^{-2}%
}+\frac{k}{k-2}}{k-1+a\left(  k-2\right)  ^{2}z^{-2}}.
\]

Letting $u=z^{-2},$ we see that%
\[%
{\displaystyle\sum\limits_{j=0}^{\infty}}
\mu_{2j}\left(  k\right)  u^{j}=-\frac{1}{2}\frac{\sqrt{1-4au}+\frac{k}{k-2}%
}{k-1+a\left(  k-2\right)  ^{2}u},
\]
and therefore%
\begin{gather*}
\sqrt{1-4au}+\frac{k}{k-2}=-2\left[  k-1+a\left(  k-2\right)  ^{2}u\right]
{\displaystyle\sum\limits_{j=0}^{\infty}}
\mu_{2j}u^{j}\\
=-%
{\displaystyle\sum\limits_{j=0}^{\infty}}
2\left(  k-1\right)  \mu_{2j}u^{j}-%
{\displaystyle\sum\limits_{j=1}^{\infty}}
2a\left(  k-2\right)  ^{2}\mu_{2\left(  j-1\right)  }u^{j}.
\end{gather*}
Since%
\[
\sqrt{1-4au}=%
{\displaystyle\sum\limits_{j=0}^{\infty}}
\binom{\frac{1}{2}}{j}\left(  -4au\right)  ^{j},
\]
we obtain%
\[
1+\frac{k}{k-2}=-2\left(  k-1\right)  \mu_{0},
\]
and%
\[
\binom{\frac{1}{2}}{j}\left(  -4a\right)  ^{j}=-2\left(  k-1\right)  \mu
_{2j}-2a\left(  k-2\right)  ^{2}\mu_{2\left(  j-1\right)  },\quad
j=1,2,\ldots.
\]

If $k=1,$ we get%
\[
\binom{\frac{1}{2}}{j}\left(  -4a\right)  ^{j}=-2a\mu_{2\left(  j-1\right)
},\quad j=1,2,\ldots,
\]
or%
\[
\mu_{2n}\left(  1\right)  =2^{2n+1}\binom{\frac{1}{2}}{n+1}\left(  -a\right)
^{n},\quad n=0,1,\ldots.
\]

If $k=2,$ we have%
\[
\mu_{2n}\left(  2\right)  =-2^{2n-1}\binom{\frac{1}{2}}{n}\left(  -a\right)
^{n},\quad n=1,2,\ldots.
\]

If $k\neq1,2,$ we set $y_{j}=\mu_{2j},$ and obtain the recurrence%
\[
y_{j+1}=-\frac{a\left(  k-2\right)  ^{2}}{k-1}y_{j}-\frac{\left(  -4a\right)
^{j+1}}{2\left(  k-1\right)  }\binom{\frac{1}{2}}{j+1},
\]
with%
\[
y_{0}=\frac{1}{2-k}.
\]
As it is well known, the general solution of the initial value problem%
\[
y_{n+1}=c_{n}y_{n}+g_{n},\quad y_{n_{0}}=y_{0},
\]
is \cite[1.2.4]{MR2128146}%
\[
y_{n}=y_{0}%
{\displaystyle\prod\limits_{j=n_{0}}^{n-1}}
c_{j}+%
{\displaystyle\sum\limits_{k=n_{0}}^{n-1}}
\left(  g_{k}%
{\displaystyle\prod\limits_{j=k+1}^{n-1}}
c_{j}\right)  .
\]
Thus,%
\[
y_{n}=\frac{1}{2-k}\left[  -\frac{a\left(  k-2\right)  ^{2}}{k-1}\right]
^{n}-%
{\displaystyle\sum\limits_{j=0}^{n-1}}
\frac{\left(  -4a\right)  ^{j+1}}{2\left(  k-1\right)  }\binom{\frac{1}{2}%
}{j+1}\left[  -\frac{a\left(  k-2\right)  ^{2}}{k-1}\right]  ^{n-j-1},
\]
or%
\[
y_{n}=-\frac{1}{2\left(  k-1\right)  }\left[  -\frac{a\left(  k-2\right)
^{2}}{k-1}\right]  ^{n}\left(  \frac{k}{k-2}+%
{\displaystyle\sum\limits_{j=0}^{n}}
\binom{\frac{1}{2}}{j}\left[  \frac{4\left(  k-1\right)  }{\left(  k-2\right)
^{2}}\right]  ^{j}\right)
\]
and the result follows.
\end{proof}

\begin{remark}
If $k=0,$ we get from (\ref{moments})%
\[
\mu_{2n}\left(  0\right)  =\frac{\left(  4a\right)  ^{n}}{2}%
{\displaystyle\sum\limits_{j=0}^{n}}
\left(  -1\right)  ^{j}\binom{\frac{1}{2}}{j},
\]
and using the identity \cite[26.3.10]{MR2723248}%
\[%
{\displaystyle\sum\limits_{j=0}^{n}}
\left(  -1\right)  ^{j}\binom{\alpha}{j}=\left(  -1\right)  ^{n}\binom
{\alpha-1}{n},
\]
we obtain%
\[
\mu_{2n}\left(  0\right)  =2^{2n-1}\binom{-\frac{1}{2}}{n}\left(  -a\right)
^{n}.
\]
This agrees with (\ref{Lc}), since
\[
\mu_{2n}\left(  0\right)  =\frac{1}{2\pi}%
{\displaystyle\int\limits_{-2\sqrt{a}}^{2\sqrt{a}}}
\frac{t^{2n}}{\sqrt{4a-t^{2}}}dt.
\]

When $k=1,$ we have from (\ref{Lc})%
\[
\mu_{2n}\left(  1\right)  =\frac{1}{2\pi a}%
{\displaystyle\int\limits_{-2\sqrt{a}}^{2\sqrt{a}}}
t^{2n}\sqrt{4a-t^{2}}dt,
\]
and therefore (\ref{mu1}) gives
\[
\frac{1}{2\pi a}%
{\displaystyle\int\limits_{-2\sqrt{a}}^{2\sqrt{a}}}
t^{2n}\sqrt{4a-t^{2}}dt=2^{2n+1}\binom{\frac{1}{2}}{n+1}\left(  -a\right)
^{n},
\]
which can be verified directly.

When $k=2,$ we can write (see Remark \ref{Remark3})%
\[
\mu_{2n}\left(  2\right)  =\frac{1}{2\pi}%
{\displaystyle\int\limits_{-2\sqrt{a}}^{2\sqrt{a}}}
t^{2n}\frac{\sqrt{4a-t^{2}}}{t^{2}}dt=\frac{1}{2\pi}%
{\displaystyle\int\limits_{-2\sqrt{a}}^{2\sqrt{a}}}
t^{2\left(  n-1\right)  }\sqrt{4a-t^{2}}dt,
\]
where $n=1,2,\ldots.$ Hence,%
\[
\mu_{2n}\left(  2\right)  =a\mu_{2\left(  n-1\right)  }\left(  1\right)
=a2^{2\left(  n-1\right)  +1}\binom{\frac{1}{2}}{n}\left(  -a\right)
^{n-1},\quad n=1,2,\ldots,
\]
in agreement with (\ref{mu2}).
\end{remark}

\section{Conclusions \label{Section4}}

We have shown that the Dickson polynomials of the $(k+1)$-th kind defined by%

\[
D_{n,k}(x;a)=%
{\displaystyle\sum\limits_{j=0}^{\left\lfloor \frac{n}{2}\right\rfloor }}
\frac{n-kj}{n-j}\binom{n-j}{j}\left(  -a\right)  ^{j}x^{n-2j}%
\]
satisfy the orthogonality relation%
\begin{gather*}
\frac{1}{2\pi}%
{\displaystyle\int\limits_{-2\sqrt{a}}^{2\sqrt{a}}}
\frac{\sqrt{4a-t^{2}}D_{n,k}\left(  t\right)  D_{m,k}\left(  t\right)
}{\left(  k-1\right)  t^{2}+a\left(  k-2\right)  ^{2}}dt\\
+\chi\left(  k\right)  \left[  \frac{1+\left(  -1\right)  ^{n+m}}{2}\right]
\frac{\left(  2-k\right)  k}{k-1}\left(  \mathrm{i}\sqrt{\frac{a}{k-1}%
}\right)  ^{n+m}=h_{n}\left(  k\right)  \delta_{n,m},
\end{gather*}
where $a>0,$ $k\in\mathbb{R},$%
\[
\chi\left(  k\right)  =\left\{
\begin{array}
[c]{c}%
0,\quad k\in\left[  0,2\right] \\
1,\quad k\in\mathbb{R}\setminus\left[  0,2\right]
\end{array}
\right.  ,
\]
and%
\[
h_{0}\left(  k\right)  =2-k,\quad h_{n}\left(  k\right)  =a^{n},\quad
n=1,2,\ldots.
\]
We hope that this work will outline some connections between finite fields and
orthogonal polynomials, and that it would be of interest to researchers in
both areas.

\begin{acknowledgement}
This paper was completed while visiting the Johannes Kepler Universit\"{a}t
Linz and supported by the strategic program "Innovatives O\"{O}-- 2010 plus"
from the Upper Austrian Government. We wish to thank Professor Peter Paule for
his generous sponsorship and our colleagues at JKU for their continuous help.

We also wish to acknowledge the hospitality of the Erwin Schr\"{o}dinger
International Institute for Mathematics and Physics (ESI), on the occasion of
the Programme on \textquotedblleft Algorithmic and Enumerative
Combinatorics\textquotedblright\ held in October-November 2017.

Finally, we also wish to express our gratitude to the anonymous referees, who
provided us with invaluable suggestions and comments that greatly improved our
first draft of the paper.
\end{acknowledgement}

\end{document}